\documentclass[secthm,seceqn,amsthm,ussrhead,12pt]{amsart}
\usepackage[utf8]{inputenc}
\usepackage[english]{babel}

\usepackage{amssymb,amsmath,amsthm,amsfonts,xcolor,enumerate,hyperref,comment,longtable,cleveref}

\usepackage{times}
\usepackage{cite}
\usepackage{pdflscape}
\usepackage{ulem}
\usepackage[mathcal]{euscript}
\usepackage{tikz}
\usepackage{hyperref}
\usepackage{cancel}
\usepackage{stmaryrd}

\usetikzlibrary{arrows}

\newtheorem{theorem}{Theorem}
\newtheorem{lemma}[theorem]{Lemma}

\theoremstyle{definition}

\theoremstyle{remark}

\newenvironment{Proof}[1][Proof.]{\begin{trivlist}
\item[\hskip \labelsep {\bfseries #1}]}{\flushright
$\Box$\end{trivlist}}

\usepackage{stmaryrd}
\usepackage{xcolor}

\newcommand{\A}{\mathbf{A}}

\begin{document}

\noindent{\large 
Degenerations of nilpotent associative commutative algebras}
\footnote{
The work was partially supported by  FAPESP 16/16445-0, 18/15712-0;
RFBR 18-31-20004; 
the President's "Program Support of Young Russian Scientists" (grant MK-2262.2019.1);
CMUP (UID/MAT/00144/2019), which is funded by FCT with national (MCTES) and European structural funds through the programs FEDER, under the partnership agreement PT2020).
The authors thank  Prof. Dr. Yury Volkov for   constructive discussions about degenerations of algebras.
} \footnote{Corresponding Author: kaygorodov.ivan@gmail.com}  

   \
   
   \ 

   {\bf      
            Ivan   Kaygorodov$^{a},$ 
            Samuel A.\ Lopes$^{b}$ \& 
            Yury Popov$^{c}$
 }
 \ 
 
\

{\tiny

$^{a}$ CMCC, Universidade Federal do ABC, Santo Andr\'e, Brasil

$^{b}$ CMUP, Faculdade de Ci\^encias, Universidade do Porto, Rua do Campo Alegre 687, 4169-007 Porto, Portugal

$^{c}$ IMECC, Universidade Estadual de Campinas, Campinas, Brasil

\ 

\

\smallskip

   E-mail addresses:

\smallskip

Ivan   Kaygorodov (kaygorodov.ivan@gmail.com) \smallskip

Samuel A.\ Lopes (slopes@fc.up.pt) \smallskip

Yury Popov (yuri.ppv@gmail.com) \smallskip

}

\

\noindent{\bf Abstract}: 
{\it 
We give a complete description of degenerations of complex
$5$-dimensional nilpotent associative commutative algebras.}

\


\ 


 \

\section*{Introduction}

There are many results related to the algebraic and geometric 
classification
of low-dimensional algebras in the varieties of Jordan, Lie, Leibniz and 
Zinbiel algebras;
for algebraic classifications  see, for example, \cite{gkks, fkkv, ikm19, gkks,  ikv18, kkk18, kpv19, ack, maz,kv16};
for geometric classifications and descriptions of degenerations see, for example, 
\cite{fkkv, ack, ale, ale2, aleis, maria, bb14, bb09, BC99, cfk19, gkks, gkk19, gkp, GRH, GRH2, ikv17, ikv18, kkk18, kpv19, kppv, kpv, kv16, kv17, S90,   gorb91, gorb93, ikm19, khud15, maz, khud13}.
Degenerations of algebras is an interesting subject, which has been studied in various papers.
In particular, there are many results concerning degenerations of algebras of small dimensions in a  variety defined by a set of identities.
One of important problems in this direction is a description of so-called rigid algebras. 
These algebras are of big interest, since the closures of their orbits under the action of the generalized linear group form irreducible components of the variety under consideration
(with respect to the Zariski topology). 
For example, rigid algebras in the varieties of
all $4$-dimensional Leibniz algebras \cite{ikv17},
all $4$-dimensional nilpotent Novikov algebras \cite{kkk18},
all $4$-dimensional nilpotent assosymmetric algebras \cite{ikm19},
all $4$-dimensional nilpotent bicommutative algebras \cite{kpv19},
all $6$-dimensional nilpotent binary Lie algebras \cite{ack},
and in some other varieties were classified.
There are fewer works in which the full information about degenerations was given for some variety of algebras.
This problem was solved 
for $2$-dimensional pre-Lie algebras  \cite{bb09},  
for $2$-dimensional terminal algebras  \cite{cfk19},
for $3$-dimensional Novikov algebras  \cite{bb14},  
for $3$-dimensional Jordan algebras  \cite{gkp},  
for $3$-dimensional Jordan superalgebras  \cite{maria},
for $3$-dimensional Leibniz and $3$-dimensional anticommutative algebras   \cite{ikv18},
for $4$-dimensional Lie algebras  \cite{BC99},
for $4$-dimensional Lie superalgebras \cite{aleis},
for $4$-dimensional Zinbiel and  $4$-dimensional nilpotent Leibniz algebras \cite{kppv},
for $5$-dimensional nilpotent Tortkara algebras  \cite{gkks},
for $6$-dimensional nilpotent Lie algebras  \cite{S90,GRH}, 
for $6$-dimensional nilpotent  Malcev algebras  \cite{kpv}, 
for $7$-dimensional $2$-step nilpotent Lie algebras \cite{ale2}, 
and for all $2$-dimensional algebras  \cite{kv16}.
Here we construct the graphs of primary degenerations for the variety of complex $5$-dimensional nilpotent
associative commutative algebras.

\newpage 
\section{Degenerations of algebras}

\subsection{Preliminaries}
Given an $n$-dimensional vector space ${\bf V}$, the set ${\rm Hom}({\bf V} \otimes {\bf V},{\bf V}) \cong {\bf V}^* \otimes {\bf V}^* \otimes {\bf V}$ 
is a vector space of dimension $n^3$. This space has a structure of the affine variety $\mathbb{C}^{n^3}.$ 
Indeed, let us fix a basis $e_1,\dots,e_n$ of ${\bf V}$. Then any $\mu\in {\rm Hom}({\bf V} \otimes {\bf V},{\bf V})$ is determined by $n^3$ structure constants $c_{i,j}^k\in\mathbb{C}$ such that
$\mu(e_i\otimes e_j)=\sum_{k=1}^nc_{i,j}^ke_k$. A subset of ${\rm Hom}({\bf V} \otimes {\bf V},{\bf V})$ is {\it Zariski-closed} if it can be defined by a set of polynomial equations in the variables $c_{i,j}^k$ ($1\le i,j,k\le n$).

Let $T$ be a set of polynomial identities.
All algebra structures on ${\bf V}$ satisfying polynomial identities from $T$ form a Zariski-closed subset of the variety ${\rm Hom}({\bf V} \otimes {\bf V},{\bf V})$. We denote this subset by $\mathbb{L}(T)$.
The general linear group ${\rm GL}({\bf V})$ acts on $\mathbb{L}(T)$ by conjugation:
$$ (g * \mu )(x\otimes y) = g\mu(g^{-1}x\otimes g^{-1}y)$$ 
for $x,y\in {\bf V}$, $\mu\in \mathbb{L}(T)\subset {\rm Hom}({\bf V} \otimes {\bf V},{\bf V})$ and $g\in {\rm GL}({\bf V})$.
Thus, $\mathbb{L}(T)$ is decomposed into ${\rm GL}({\bf V})$-orbits that correspond to the isomorphism classes of algebras. 
Let $O(\mu)$ denote the ${\rm GL}({\bf V})$-orbit of $\mu\in\mathbb{L}(T)$ and $\overline{O(\mu)}$ its Zariski closure.

Let ${\bf A}$ and ${\bf B}$ be two $n$-dimensional algebras satisfying identities from $T$ and $\mu,\lambda \in \mathbb{L}(T)$ represent ${\bf A}$ and ${\bf B}$ respectively.
We say that ${\bf A}$ {\it degenerates to} ${\bf B}$ and write ${\bf A}\to {\bf B}$ if $\lambda\in\overline{O(\mu)}$.
Note that in this case we have $\overline{O(\lambda)}\subset\overline{O(\mu)}$. Hence, the definition of a degeneration does not depend on the choice of $\mu$ and $\lambda$. If ${\bf A}\not\cong {\bf B}$, then the assertion ${\bf A}\to {\bf B}$ 
is called a {\it proper degeneration}. We write ${\bf A}\not\to {\bf B}$ if $\lambda\not\in\overline{O(\mu)}$.

Let ${\bf A}$ be represented by $\mu\in\mathbb{L}(T)$. Then  ${\bf A}$ is  {\it rigid} in $\mathbb{L}(T)$ if $O(\mu)$ is an open subset of $\mathbb{L}(T)$.
Recall that a subset of a variety is called {\it irreducible} if it cannot be represented as a union of two non-trivial closed subsets. A maximal irreducible closed subset of a variety is called an {\it irreducible component}.
It is well known that any affine variety can be represented as a finite union of its irreducible components in a unique way.
The algebra ${\bf A}$ is rigid in $\mathbb{L}(T)$ if and only if $\overline{O(\mu)}$ is an irreducible component of $\mathbb{L}(T)$.

In the present work we use the methods applied to Lie algebras in \cite{BC99,GRH,GRH2,S90}.
First of all, if ${\bf A}\to {\bf B}$ and ${\bf A}\not\cong {\bf B}$, then $\dim \mathfrak{Der}({\bf A})<\dim \mathfrak{Der}({\bf B})$, where $\mathfrak{Der}({\bf A})$ is the Lie algebra of derivations of ${\bf A}$. We will compute the dimensions of algebras of derivations and will check the assertion ${\bf A}\to {\bf B}$ only for such ${\bf A}$ and ${\bf B}$ that $\dim \mathfrak{Der}({\bf A})<\dim \mathfrak{Der}({\bf B})$. Secondly, if ${\bf A}\to {\bf C}$ and ${\bf C}\to {\bf B}$ then ${\bf A}\to{\bf  B}$. If there is no ${\bf C}$ such that ${\bf A}\to {\bf C}$ and ${\bf C}\to {\bf B}$ are proper degenerations, then the assertion ${\bf A}\to {\bf B}$ is called a {\it primary degeneration}. If $\dim \mathfrak{Der}({\bf A})<\dim \mathfrak{Der}({\bf B})$ and there are no ${\bf C}$ and ${\bf D}$ such that ${\bf C}\to {\bf A}$, ${\bf B}\to {\bf D}$, ${\bf C}\not\to {\bf D}$ and one of the assertions ${\bf C}\to {\bf A}$ and ${\bf B}\to {\bf D}$ is a proper degeneration,  then the assertion ${\bf A} \not\to {\bf B}$ is called a {\it primary non-degeneration}. It suffices to prove only primary degenerations and non-degenerations to describe degenerations in the variety under consideration. It is easy to see that any algebra degenerates to the algebra with zero multiplication. From now on we use this fact without mentioning it.

To prove primary degenerations, we will construct families of matrices parametrized by $t$. Namely, let ${\bf A}$ and ${\bf B}$ be two algebras represented by the structures $\mu$ and $\lambda$ from $\mathbb{L}(T)$ respectively. Let $e_1,\dots, e_n$ be a basis of ${\bf V}$ and $c_{i,j}^k$ ($1\le i,j,k\le n$) be the structure constants of $\lambda$ in this basis. If there exist $a_i^j(t)\in\mathbb{C}$ ($1\le i,j\le n$, $t\in\mathbb{C}^*$) such that $E_i^t=\sum_{j=1}^na_i^j(t)e_j$ ($1\le i\le n$) form a basis of ${\bf V}$ for any $t\in\mathbb{C}^*$, and the structure constants $c_{i,j}^k(t)$ of $\mu$ in the basis $E_1^t,\dots, E_n^t$ satisfy $\lim\limits_{t\to 0}c_{i,j}^k(t)=c_{i,j}^k$, then ${\bf A}\to {\bf B}$. In this case  $E_1^t,\dots, E_n^t$ is called a {\it parametric basis} for ${\bf A}\to {\bf B}$.

To prove primary non-degenerations we will use the following lemma (see \cite{GRH}).

\begin{lemma}\label{main}
Let $\mathcal{B}$ be a Borel subgroup of ${\rm GL}({\bf V})$ and $\mathcal{R}\subset \mathbb{L}(T)$ be a $\mathcal{B}$-stable closed subset.
If ${\bf A} \to {\bf B}$ and ${\bf A}$ can be represented by $\mu\in\mathcal{R}$ then there is $\lambda\in \mathcal{R}$ that represents ${\bf B}$.
\end{lemma}
 
Each time when we will need to prove some primary non-degeneration $\mu\not\to\lambda$, we will define $\mathcal{R}$ by a set of polynomial equations in structure constants $c_{ij}^k$ in such a way that the structure constants of $\mu$ in the basis $e_1,\dots,e_n$ satisfy these equations. We will omit everywhere the verification of the fact that $\mathcal{R}$ is stable under the action of the subgroup of lower triangular matrices and of the fact that $\lambda\not\in\mathcal{R}$ for any choice of a basis of ${\bf V}.$ 
To simplify our equations, we will use the notation $A_i=\langle e_i,\dots,e_n\rangle,\ i=1,\ldots,n$ and write simply $A_pA_q\subset A_r$ instead of $c_{ij}^k=0$ ($i\geq p$, $j\geq q$, $k\geq r$).

If the number of orbits under the action of ${\rm GL}({\bf V})$ on  $\mathbb{L}(T)$ is finite, then the graph of primary degenerations gives the whole picture. In particular, the description of rigid algebras and irreducible components can be easily obtained. 

\subsection{Degenerations of $5$-dimensional nilpotent associative commutative algebras}
The algebraic classification of $5$-dimensional nilpotent associative commutative algebras was given in \cite{maz}.
Also, in the same paper, it was proved that the variety of all $5$-dimensional nilpotent associative commutative algebras 
has only one irreducible component.
The main result of the present section is the following theorem.

\begin{theorem}
The graph of all degenerations in the variety of $5$-dimensional nilpotent  associative commutative algebras is given in Figure B (see, Appendix).
\end{theorem}

\begin{Proof} Tables C, D 
presented in Appendix give the proofs for all primary degenerations and non-degenerations.
\end{Proof}

\section{Appendix.}

\begin{center}

$$	\begin{array}{|l|c|llllll|}
			\hline 
			\multicolumn{8}{|c|}{\textrm{{\bf Table A. $5$-dimensional nilpotent associative commutative  algebras}}}  \\
			\hline
			\hline 

			\A & \mathfrak{Der} \ \A & 			\multicolumn{6}{|c|}{\textrm{Multiplication table}}  \\

			\hline
			
			\A_{01} & 5 & e_{1}^{2}=e_{2},& e_{2}^{2}=e_4, &e_{1}e_{3}=e_{4},& e_{1}e_{2}=e_{3},& e_{1}e_{4}=e_5, &e_{2}e_{3}=e_{5}\\
            \hline
            \A_{02} & 6 & e_{1}^{2}=e_{3},&  e_{2}^{2}= e_5, & e_3^2 = e_5,& e_{1}e_{3}=e_{4}, & e_1e_4 = e_5 &\\
            \hline
            \A_{03} & 6 & e_{1}^{2}=e_{3},& e_{2}^{2}=e_{4},& e_{1}e_{3}=e_5, &e_{2}e_{4}=e_{5} &&\\
            \hline
            \A_{04} & 7 & e_{1}^{2}=e_{3},& e_{1}e_{2}=e_{4},& e_{1}e_{4}=e_5,& e_{2}e_{3}=e_{5} &&\\
            \hline
            \A_{05} & 7 & e_{1}^{2}=e_{2},& e_{2}^{2}=e_{4},& e_{1}e_{2}=e_{3},& e_{1}e_{3}=e_{4} &&\\
            \hline
            \A_{06} & 7 & e_{1}^{2}=e_{2},& e_{1}e_{2}=e_{3},& e_4^2 = e_5 &&&\\
            \hline
            \A_{07} & 7 &e_{1}e_{3}=e_{4},& e_{2}e_{3}=e_{5},& e_{1}e_{2}=e_{4}+e_{5} &&&\\
            \hline
            \A_{08} & 8 & e_{1}^{2}=e_{3},&  e_{2}^{2}=e_4, & e_{1}e_{3}=e_{4},& e_{1}e_{2}=e_{5} &&\\
            \hline
            \A_{09} & 8 & e_{1}e_{3}=e_{5},&  e_{1}e_{2}=e_{4},& e_{2}e_{3}=-e_{5} &&&\\
            \hline
            \A_{10} & 9 & e_{1}^{2}=e_{3},& e_{1}e_{3}=e_{4},& e_{1}e_{2}=e_{5} &&&\\
            \hline
            \A_{11} & 9 & e_{1}^{2}=e_4, &e_{2}e_{3}=e_{4},& e_{1}e_{3}=e_{5} &&&\\
            \hline
            \A_{12} & 11 & e_{1}e_{2}=e_{4},&  e_{1}e_{3}=e_{5} &&&&\\
            \hline
            \A_{13} & 8 & e_{3}^{2}=e_{4},& e_{1}e_{2}=e_5, & e_{3}e_{4}=e_{5} &&&\\
            \hline
            \A_{14} & 9 & e_{1}^{2}=e_{3},& e_{2}^{2}=e_{4},& e_{1}e_{3}=e_{4} &&&\\
            \hline
            \A_{15} & 9 & e_1e_2 = e_3,& e_4^2 = e_5 &&&&\\
            \hline
            \A_{16} & 10 & e_{1}^{2}=e_{3},& e_{2}^{2}=e_{5},& e_{1}e_{2}=e_{4} &&&\\
            \hline
            \A_{17} & 10 & e_{1}^{2}=e_{4},&  e_{3}^{2}=e_5, &e_{1}e_{2}=e_{5} &&&\\
            \hline
            \A_{18} & 11 & e_{1}^{2}=e_{2},& e_{1}e_{2}=e_{3} &&&&\\
            \hline
            \A_{19} & 11 & e_{1}^{2}=e_{3},& e_2^2 = e_4 &&&&\\
            \hline
            \A_{20} & 12 & e_1^2 = e_3,& e_1e_2 = e_4 &&&&\\
            \hline
            \A_{21} & 11 & e_2e_3 =e_5, & e_1e_4 = e_5 &&&&\\
            \hline
            \A_{22} & 12 & e_1^2 = e_4,& e_2e_3 = e_4 &&&&\\
            \hline
            \A_{23} & 14 & e_1e_2 = e_3 &&&&&\\
            \hline
            \A_{24} & 17 & e_1^2 = e_2 &&&&&\\
            \hline
%            \mathbb{C}^5 & 25 &  &&&&\\
 %           \hline
			\end{array}$$ 

\end{center}

\newpage
\begin{center}
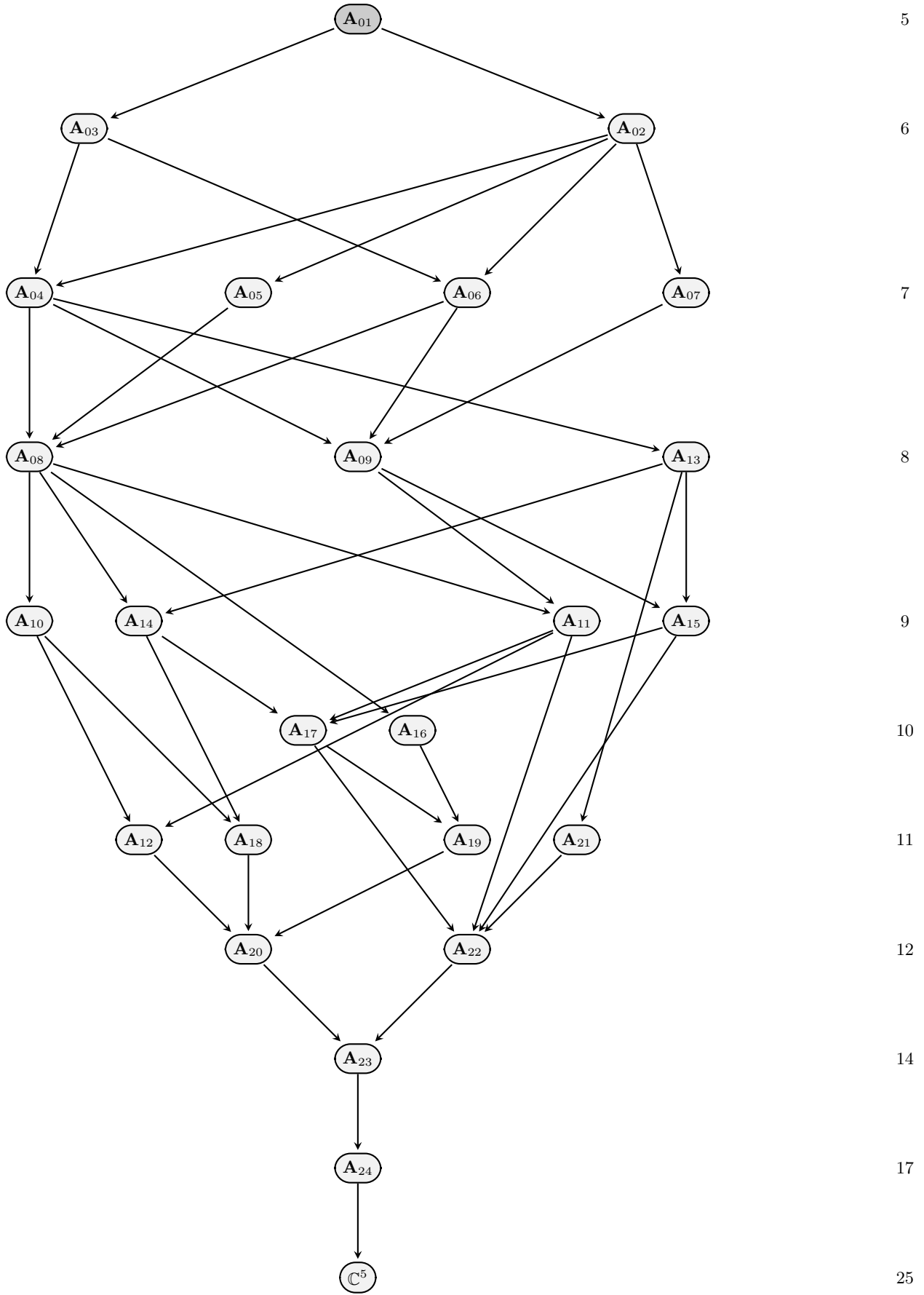
	{\bf Figure B.  The graph of degenerations of $5$-dimensional nilpotent associative commutative algebras}
	
	\
	
	\begin{tikzpicture}[->,>=stealth,shorten >=0.05cm,auto,node distance=1.3cm,
	thick,main node/.style={rectangle,draw,fill=gray!10,rounded corners=1.5ex,font=\sffamily \scriptsize \bfseries },rigid node/.style={rectangle,draw,fill=black!20,rounded corners=1.5ex,font=\sffamily \scriptsize \bfseries },style={draw,font=\sffamily \scriptsize \bfseries }]

	\node (5) at (0,24) {$5$};
	\node (6) at (0,22) {$6$};
	\node (7) at (0,19) {$7$};
	\node (8) at (0,16) {$8$};
	\node (9) at (0,13) {$9$};
	\node (10) at (0,11) {$10$};
	\node (11)  at (0,9) {$11$};
	\node (12)  at (0,7) {$12$};
	\node (14)  at (0,5) {$14$};
	\node (17)  at (0,3) {$17$};
	\node (25)  at (0,1) {$25$};

	\node[rigid node] (a01) at (-10,24) {$\A_{01}$};
	
	\node[main node] (a02) at (-5,22) {$\A_{02}$};
	\node[main node] (a03) at (-15,22) {$\A_{03}$};
	
	\node[main node] (a04) at (-16,19) {$\A_{04}$};
    \node[main node] (a05) at (-12,19) {$\A_{05}$};
	\node[main node] (a06) at (-8,19) {$\A_{06}$};
	\node[main node] (a07) at (-4,19) {$\A_{07}$};

	\node[main node] (a08) at (-16,16) {$\A_{08}$};
	\node[main node] (a09) at (-10,16) {$\A_{09}$};
	\node[main node] (a13) at (-4,16) {$\A_{13}$};

	\node[main node] (a10) at (-16,13) {$\A_{10}$};
	\node[main node] (a14) at (-14,13) {$\A_{14}$};
	\node[main node] (a11) at (-6,13) {$\A_{11}$};
	\node[main node] (a15) at (-4,13) {$\A_{15}$};
	
	\node[main node] (a16) at (-9,11) {$\A_{16}$};
	\node[main node] (a17) at (-11,11) {$\A_{17}$};
	
	\node[main node] (a12) at (-14,9) {$\A_{12}$};
	\node[main node] (a18) at (-12,9) {$\A_{18}$};
	\node[main node] (a19) at (-8,9) {$\A_{19}$};
	\node[main node] (a21) at (-6,9) {$\A_{21}$};
	
	\node[main node] (a20) at (-12,7) {$\A_{20}$};
	\node[main node] (a22) at (-8,7) {$\A_{22}$};
	
	\node[main node] (a23) at (-10,5) {$\A_{23}$};
	
	\node[main node] (a24) at (-10,3) {$\A_{24}$};
	
	\node[main node] (CC)  at (-10,1) {$\mathbb{C}^5$};
	
	\path[every node/.style={font=\sffamily\small}]
   
 (a01) edge   node[left] {} (a02)
 (a01) edge   node[left] {} (a03)
 
 (a02) edge   node[left] {} (a04)
 (a02) edge   node[left] {} (a05)
 (a02) edge   node[left] {} (a06)
 (a02) edge   node[left] {} (a07)
 
 (a03) edge   node[left] {} (a04)
 (a03) edge   node[left] {} (a06)
 
 (a04) edge   node[left] {} (a08)
 (a04) edge   node[left] {} (a09)
 (a04) edge   node[left] {} (a13)
 
% (a05) edge   node[left] {} (a09)
 (a05) edge   node[left] {} (a08)
 
 (a06) edge   node[left] {} (a08)
 (a06) edge   node[left] {} (a09)
% (a06) edge   node[left] {} (a14)
% (a06) edge   node[left] {} (a16)

 (a07) edge   node[left] {} (a09)
 %(a07) edge   node[left] {} (a11)

 (a08) edge   node[left] {} (a10)
 (a08) edge   node[left] {} (a11)
 (a08) edge   node[left] {} (a14)
 (a08) edge   node[left] {} (a16)
 
 (a09) edge   node[left] {} (a11)
 %(a09) edge   node[left] {} (a12)
 (a09) edge   node[left] {} (a15)
 %(a09) edge   node[left] {} (a17)
   
 (a10) edge   node[left] {} (a12)
 (a10) edge   node[left] {} (a18)

 (a11) edge   node[left] {} (a12)
 (a11) edge   node[left] {} (a17)
 (a11) edge   node[left] {} (a22)

 (a12) edge   node[left] {} (a20)

 (a13) edge   node[left] {} (a15)
 (a13) edge   node[left] {} (a14)
 (a13) edge   node[left] {} (a21)

 (a14) edge   node[left] {} (a17)
 (a14) edge   node[left] {} (a18)
 %(a14) edge   node[left] {} (a22)
 
 (a15) edge   node[left] {} (a17)
 (a15) edge   node[left] {} (a22)

 (a16) edge   node[left] {} (a19)

 (a17) edge   node[left] {} (a19)
 (a17) edge   node[left] {} (a22)

 (a18) edge   node[left] {} (a20)
 
 (a19) edge   node[left] {} (a20)
 
 (a20) edge   node[left] {} (a23)
 
 (a21) edge   node[left] {} (a22)

 (a22) edge   node[left] {} (a23)
 
 (a23) edge   node[left] {} (a24)
 
 (a24) edge   node[left] {} (CC)
    
    ;
   
	\end{tikzpicture}
	
\end{center}

\newpage

{\tiny
$$\begin{array}{|lll|}
\hline

\multicolumn{3}{|c|}{\textrm{{\bf Table C. Degenerations of $5$-dimensional nilpotent associative commutative  algebras}}}  \\
\hline
%\multicolumn{3}{|c|}{\textrm{Degeneration}}  & \multicolumn{5}{|c|}{\textrm{Parametric basis}} \\
%\hline
%hline

\begin{array}{l|}\A_{01} \to \A_{02} \\
\hline\end{array} & 
E_{1}^t= t e_1&
E_{2}^t= -i (t^{-1} e_3 - t^{-4} e_4 + t^{-7} e_5 - t^2 e_2)\\
E_{3}^t= t^2 e_2&
E_{4}^t= e_4 - t^{-3} e_5 &
E_{5}^t= t e_5\\
\hline

\begin{array}{l|}\A_{01} \to \A_{03} \\
\hline\end{array}   & 
\begin{array}{l} 
E_{1}^{t}=t e_1 + \frac{2t-1}{2-3t}e_2 + \\ \frac{1 - 5 t + 5 t^2}{2t(2-3t)^2}e_3 + \frac{-3 + 9 t - 7 t^2}{4t(-2+3t)^3}e_4\end{array}& 
\begin{array}{l}E_{2}^{t}=t e_1 + \frac{1-t}{2-3t}e_2 -\\ \frac{-1 + t + t^2}{2t(2-3t)^2}e_3 + \frac{-3 + 9 t - 7 t^2}{4t(-2+3t)^3}e_4 \end{array}\\ 
\begin{array}{l} E_{3}^{t}=-te_3 +\frac{1-3t}{2-3t}e_4 + \\ \frac{-1 + 7 t - 9 t^2}{2t(2-3t)^2}e_5 \end{array} & 
E_{4}^{t}=te_3 +\frac{1}{2-3t}e_4 + \frac{1 + t - 3t^2}{2t(2-3t)^2}e_5& 
E_{5}^{t}=t e_5\\
\hline  

%\A_{01} \to \A_{07}  & \begin{array}{l} 
%E_{1}^{t}= t e_1 - \frac{1}{t} e_3 + \frac{1 - t}{2 t^2} e_4\end{array}& 
%\begin{array}{l} E_{2}^{t}= t e_2 - e_3 \end{array} & 
%\begin{array}{l} E_{3}^{t}= e_3 - \frac{2}{t} e_4 + \frac{1 - t}{t^2} e_5\end{array}  &
%E_{4}^{t}= t e_4 - 2 e_5& 
%E_{5}^{t}= t e_5\\
%\hline

\begin{array}{l|}\A_{02} \to \A_{04}   \\
\hline\end{array} & 

E_{1}^t= t e_1 &
E_{2}^t= i e_2 + e_3\\
E_{3}^t= t^2 e_3 &
E_{4}^t= t e_4 &
E_{5}^t= t^2 e_5 \\

%E_{1}^{t}=t e_1& 
%E_{2}^{t}=-i(t e_4-t^2 e_3) \\ 
%E_{3}^{t}= t^2 e_3  &
%E_{4}^{t}=i e_2 + e_3 & 
%E_{5}^{t}=t^2 e_5\\
\hline

\begin{array}{l|}\A_{02} \to \A_{05}  \\
\hline\end{array}  &  
E_{1}^{t}=e_1& 
E_{2}^{t}=e_3  \\ 
E_{3}^{t}=e_4  &
E_{4}^{t}=e_5& 
E_{5}^{t}=te_2\\
\hline

\begin{array}{l|}\A_{02} \to \A_{06}  \\
\hline\end{array}  & 
E_{1}^{t}=e_1& 
E_{2}^{t}= e_3 \\
E_{3}^{t}= e_4& 
E_{4}^{t}=\sqrt{\frac{-1 - t^3}{t}} e_2 + t e_3  &  
E_{5}^{t}=-\frac{1}{t} e_5\\
\hline

\begin{array}{l|}\A_{02} \to \A_{07}  \\
\hline\end{array}  &  
E_{1}^{t}=te_1 + \frac{1}{3}e_3 & 
E_{2}^{t}=ie_2 + e_3 + \frac{5}{3t}e_4  \\
E_{3}^{t}= -ie_2 + e_3 - \frac{1}{3t}e_4 + \frac{2}{9t^2}e_5  &
E_{4}^{t}= te_4& 
E_{5}^{t}=2e_5\\
\hline

\begin{array}{l|} \A_{03} \to \A_{04}  \\
\hline\end{array}  &  
E_{1}^{t}=  t e_1 - t e_2& 
 E_{2}^{t}= t^2 e_2  \\  
 E_{3}^{t}= t^2 e_3 + t^2 e_4 &
E_{4}^{t}= -t^3 e_4 &
E_{5}^{t}= t^4 e_5\\
\hline

\begin{array}{l|} \A_{03} \to \A_{06} \\
\hline\end{array}   &  
E_{1}^{t}=  e_1 & 
E_{2}^{t}=  e_3\\
E_{3}^{t}=  e_5 &
E_{4}^{t}= t e_2 & 
E_{5}^{t}= t^2 e_4\\
\hline

\begin{array}{l|} \A_{04} \to \A_{08}   \\
\hline\end{array} &  
E_{1}^{t}= e_1 + \frac{1}{3t}e_2 - \frac{1}{6t^2}e_3 + \frac{1}{18t^3}e_4& 
E_{2}^{t}= e_2 + \frac{1}{2t}e_3   \\
E_{3}^{t}= e_3 + \frac{2}{3t}e_4  &
E_{4}^{t}= \frac{1}{t}e_5 &
E_{5}^{t}= e_4\\
\hline

\begin{array}{l|} \A_{04} \to \A_{09}  \\
\hline\end{array}  &  
E_{1}^{t}= te_1 - \frac{1}{3}e_2 & 
E_{2}^{t}= e_2  \\
E_{3}^{t}= te_3 - \frac{2}{3}e_4  &
E_{4}^{t}= te_4 &
E_{5}^{t}= -te_5\\
\hline

\begin{array}{l|} \A_{04} \to \A_{13}  \\
\hline\end{array}  &  
E_{1}^{t}=  e_2& 
E_{2}^{t}= te_3 - \frac{1}{3t}e_4  \\ 
E_{3}^{t}= te_1 + \frac{1}{3t}e_2   &
E_{4}^{t}= e_4 &
E_{5}^{t}= te_5\\
\hline

\begin{array}{l|} \A_{05} \to \A_{08}  \\
\hline\end{array}  &  
E_{1}^{t}= t e_1& 
 E_{2}^{t}= t e_2 - e_3 - e_5  \\ 
 E_{3}^{t}=t e_3 + t e_5  &
E_{4}^{t}= -t^3 e_5 &
E_{5}^{t}= -t e_4 - t^2 e_5\\
\hline

%\begin{array}{l|} \A_{05} \to \A_{09}  \\
%\hline\end{array}  &  
%E_{1}^{t}= t e_1& 
%% E_{2}^{t}= -t e_1+t e_2 - t^2  e_3+  (1 + t^3)  e_4 + e_5 \\
 % E_{3}^{t}= e_3  -\frac{1 + t^3}{t^2} e_4 - \frac{1}{t^2} e_5  &
%E_{4}^{t}= e_4 +  e_5 &
%E_{5}^{t}= -t e_5\\
%\hline

\begin{array}{l|} \A_{06} \to \A_{08}  \\
\hline\end{array} &   
E_{1}^t= e_1 + e_4 &
E_{2}^t= t e_4 - t^{-1} e_2 \\
E_{3}^t= e_2 +  e_5 &
E_{4}^t= e_3&
E_{5}^t= t e_5 - t^{-1} e_3 \\
\hline

\begin{array}{l|} \A_{06} \to \A_{09}  \\
\hline\end{array} &   
E_{1}^t= t e_1 + t e_4 &
E_{2}^t= -(t e_1 + t e_4) + 2 t e_4 + t   e_2 \\
E_{3}^t= e_2 +  e_5 &
E_{4}^t= 2 t^2 e_5 + t^2 e_3 &
E_{5}^t= t e_3\\
\hline

%\begin{array}{l|} \A_{06} \to \A_{14}  &   
%E_{1}^{t}=  e_1& 
%  E_{2}^{t}= t e_2   & 
%  E_{3}^{t}= e_3  &
%E_{4}^{t}= e_5 &
%E_{5}^{t}= e_4 - \frac{1}{t^2} e_5\\
%\hline

%\begin{array}{l|} \A_{06} \to \A_{16}  &   
%E_{1}^{t}= t e_1 + e_4& 
%E_{2}^{t}= t e_4 +  e_2 & 
%  E_{3}^{t}= t^2 e_2 + e_5 &
%E_{4}^{t}= t e_5 + t e_3 &
%E_{5}^{t}= -t^2 e_3 \\
%\hline

\begin{array}{l|} \A_{07} \to \A_{09} \\
\hline\end{array}  &   
E_{1}^t= -e_1 - e_2 + e_3 &
E_{2}^t= e_2 \\
E_{3}^t= t e_3 &
E_{4}^t= - e_4  &
E_{5}^t= -t e_5\\
\hline

%\begin{array}{l|} \A_{07} \to \A_{11} \\
%\hline\end{array}  &   
%E_{1}^{t}= X e_1 + \frac{1}{2}(-t + X) e_2 - \frac{1}{2}(t + X)e_3, 
%& 
%  E_{2}^{t}= t e_2 - te_3  \\ 
% X=\sqrt{4 t + t^2 - 4 t^3}, \  E_{3}^{t}= 2 e_3   &
%E_{4}^{t}= 2t e_5 &
%%E_{5}^{t}= 2X e_4 + (-t+2t^2 + X)e_5\\
%\hline

\begin{array}{l|} \A_{08} \to \A_{10} \\
\hline\end{array}  &   
E_{1}^{t}=te_1& 
  E_{2}^{t}=t^2e_2  \\ 
  E_{3}^{t}=t^2e_3  &
E_{4}^{t}=t^3e_4 &
E_{5}^{t}=t^3e_5\\
\hline

\begin{array}{l|} \A_{08} \to \A_{11} \\
\hline\end{array}  &   
E_{1}^{t}=te_2& 
  E_{2}^{t}=te_3  \\ 
  E_{3}^{t}=te_1  &
E_{4}^{t}=t^2e_4 &
E_{5}^{t}=t^2e_5\\
\hline

\begin{array}{l|} \A_{08} \to \A_{14} \\
\hline\end{array}  &   
E_{1}^{t}=e_1& 
  E_{2}^{t}=e_2  \\
  E_{3}^{t}=e_3  &
E_{4}^{t}=e_4 &
E_{5}^{t}=\frac{1}{t}e_5\\
\hline

\begin{array}{l|} \A_{08} \to \A_{16} \\
\hline\end{array}  &   
E_{1}^{t}=te_1& 
  E_{2}^{t}=te_2 + t^{-1} e_3  \\
  E_{3}^{t}=t^2e_3  &
E_{4}^{t}=e_4 + t^2e_5 &
E_{5}^{t}=-t^4e_5\\
\hline

\begin{array}{l|} \A_{09} \to \A_{11} \\
\hline\end{array}  &   
E_{1}^{t}=te_1 + \frac{1}{2t}e_3& 
  E_{2}^{t}=-e_3  \\ 
  E_{3}^{t}=e_2  &
E_{4}^{t}=e_5 &
E_{5}^{t}=te_4 - \frac{1}{2t}e_5\\
\hline

\begin{array}{l|} \A_{09} \to \A_{12} \\
\hline\end{array}  &   
E_{1}^{t}=e_1& 
  E_{2}^{t}=te_2  \\ 
  E_{3}^{t}=te_3  &
E_{4}^{t}=te_4 &
E_{5}^{t}=te_5\\
\hline

\begin{array}{l|} \A_{09} \to \A_{15} \\
\hline\end{array}  &   
E_{1}^{t}=e_3 &
  E_{2}^{t}=t^2 e_2  \\ 
  E_{3}^{t}=-t^2e_5  &
E_{4}^{t}=t  e_1+t e_2& 
E_{5}^{t}=2t^2e_4\\
\hline

%\begin{array}{l|} \A_{09} \to \A_{17}\\
%\hline\end{array}   &   
%E_{1}^{t}= t e_1 &
%  E_{2}^{t}= t e_3 + \frac{2}{t} e_4 \\
%  E_{3}^{t}= e_2 + e_3  & 
%E_{4}^{t}= t^2 e_4& 
%E_{5}^{t}= 2 e_5\\
%\hline

\begin{array}{l|} \A_{10} \to \A_{12} \\
\hline\end{array}  &   
E_{1}^{t}=te_1& 
  E_{2}^{t}=-t^{-1}e_2 + e_3 \\
  E_{3}^{t}=te_3  &
E_{4}^{t}=te_4-e_5 &
E_{5}^{t}=t e_5\\
\hline

\begin{array}{l|} \A_{10} \to \A_{18}\\
\hline\end{array}   &   
E_{1}^t= t e_1&
E_{2}^t= t^2  e_3 \\
E_{3}^t= t^3  e_4 &
E_{4}^t= e_2&
E_{5}^t= e_5\\
\hline

\begin{array}{l|} \A_{11} \to \A_{12} \\
\hline\end{array}  &   
E_{1}^{t}=e_3& 
  E_{2}^{t}=t e_2  \\
  E_{3}^{t}=te_1  & 
E_{4}^{t}=te_4 &
E_{5}^{t}=te_5\\
\hline

\begin{array}{l|} \A_{11} \to \A_{17} \\
\hline\end{array}  &   
E_{1}^t= t e_1 + \frac{1 + t^3}{4} e_2 - 2t^{-1} e_3 &
E_{2}^t= t e_3  \\
E_{3}^t= \frac{t^4}{8} e_2 +  e_3 &
E_{4}^t= -t^{-1} e_4 - 4  e_5 &
E_{5}^t= -t^5 e_5\\
\hline

\begin{array}{l|} \A_{11} \to \A_{22} \\
\hline\end{array}  &   
E_{1}^{t}=te_1& 
  E_{2}^{t}=te_2  \\
  E_{3}^{t}=t e_3  &
E_{4}^{t}=t^2e_4 &
E_{5}^{t}= e_5\\
\hline

\begin{array}{l|} \A_{12} \to \A_{20} \\
\hline\end{array}  &   
E_{1}^{t}=te_1+\frac{1}{2} e_2+te_3& 
  E_{2}^{t}=te_2+\frac{1}{t} e_3 \\
  E_{3}^{t}=te_4  &
E_{4}^{t}=e_5 &
E_{5}^{t}=e_3\\
\hline

\begin{array}{l|} \A_{13} \to \A_{14}\\
\hline\end{array}  &   
E_{1}^{t}=e_3& 
  E_{2}^{t}=te_2+\frac{1}{2t}e_1 \\
  E_{3}^{t}=e_4  &
E_{4}^{t}=e_5 &
E_{5}^{t}=e_1
\\
\hline

\end{array}$$
}

{\tiny 
$$\begin{array}{|llllll|}
\hline
\begin{array}{l|} \A_{13} \to \A_{15} \\
\hline\end{array}  &   
E_{1}^{t}=e_1& 
  E_{2}^{t}=e_2  &
  E_{3}^{t}=e_5  &
E_{4}^{t}=te_3 &
E_{5}^{t}=t^2e_4\\
\hline

\begin{array}{l|} \A_{13} \to \A_{21} \\
\hline\end{array}  &   
E_{1}^{t}=t e_3& 
  E_{2}^{t}=t e_2 &
  E_{3}^{t}=te_1  &
E_{4}^{t}=te_4 &
E_{5}^{t}=t^2e_5\\
\hline

\begin{array}{l|} \A_{14} \to \A_{17} \\
\hline\end{array}  &   
E_{1}^{t}=te_1& 
  E_{2}^{t}=te_3-e_5 & 
  E_{3}^{t}=te_2  &
E_{4}^{t}=te_5 &
E_{5}^{t}=t^2e_4\\
\hline

\begin{array}{l|} \A_{14} \to \A_{18}  \\
\hline\end{array} &   
E_{1}^{t}=te_1& 
  E_{2}^{t}=t^2e_3  &
  E_{3}^{t}=t^3e_4  &
E_{4}^{t}=t^3e_2 &
E_{5}^{t}=e_5\\
\hline

%\begin{array}{l|} \A_{14} \to \A_{22} \\
%\hline\end{array}  &   
%E_{1}^{t}=t e_2 & 
%%  E_{2}^{t}=t e_1 &
 % E_{3}^{t}=t e_3  &
%E_{4}^{t}=t^2 e_4 &
%%E_{5}^{t}=e_5\\
%\hline

\begin{array}{l|} \A_{15} \to \A_{17} \\
\hline\end{array}  &   
E_{1}^{t}=e_4 + e_1 + e_2,& 
  E_{2}^{t}=-2t^2e_2 &
  E_{3}^{t}=te_1-te_2  &
E_{4}^{t}=e_5 + 2e_3 &
E_{5}^{t}=-2t^2e_3\\
\hline

\begin{array}{l|} \A_{15} \to \A_{22} \\
\hline\end{array}  &   
E_{1}^{t}=te_4& 
  E_{2}^{t}=e_2 &
  E_{3}^{t}=te_1  &
E_{4}^{t}=te_3 &
E_{5}^{t}=t e_5-e_3\\
\hline

\begin{array}{l|} \A_{16} \to \A_{19} \\
\hline\end{array}  &   
E_{1}^{t}=te_1& 
  E_{2}^{t}=te_2  &
  E_{3}^{t}=t^2e_3   &
E_{4}^{t}=t^2e_5 &
E_{5}^{t}=e_4\\
\hline

\begin{array}{l|} \A_{17} \to \A_{19} \\
\hline\end{array}  &   
E_{1}^t= t e_1 - t^2 e_3 &
E_{2}^t= t e_2 +  e_3&
E_{3}^t= t^2 e_4 + t^4 e_5 &
E_{4}^t= e_5&
E_{5}^t= t e_3\\
\hline

\begin{array}{l|} \A_{17} \to \A_{22} \\
\hline\end{array}  &   
E_{1}^t= -2te_1 - \frac{1}{t} e_2 - e_3 &
E_{2}^t= t e_1 +  e_3&
E_{3}^t= -2te_1 - \frac{1}{t} e_2 &
E_{4}^t= e_5&
E_{5}^t= e_4 + \frac{1}{t^2}e_5\\
\hline

\begin{array}{l|} \A_{18} \to \A_{20} \\
\hline\end{array}  &   
E_{1}^{t}=te_1& 
  E_{2}^{t}=te_2-e_3+\frac{1}{t} e_4 &
  E_{3}^{t}=t e_3-e_4  &
E_{4}^{t}=t e_4 &
E_{5}^{t}=e_5\\
\hline

\begin{array}{l|} \A_{19} \to \A_{20}\\
\hline\end{array}   &   
E_{1}^{t}=te_1+e_2& 
  E_{2}^{t}=te_2  &
  E_{3}^{t}=t^2 e_3+e_4  &
E_{4}^{t}=t e_4 &
E_{5}^{t}=e_5\\
\hline

\begin{array}{l|} \A_{20} \to \A_{23}  \\
\hline\end{array} &   
E_{1}^{t}=e_1 & 
  E_{2}^{t}=e_2  &
  E_{3}^{t}=e_4  &
E_{4}^{t}=\frac{1}{t}e_3 &
E_{5}^{t}=e_5\\
\hline

\begin{array}{l|} \A_{21} \to \A_{22}  \\
\hline\end{array} &   
E_{1}^{t}=te_1+ \frac{1}{2} te_4& 
  E_{2}^{t}=te_2  &
  E_{3}^{t}=t e_3  &
E_{4}^{t}=t^2 e_5 &
E_{5}^{t}=t^2 e_4 \\
\hline

\begin{array}{l|} \A_{22} \to \A_{23} \\
\hline\end{array}  &   
E_{1}^{t}=e_2& 
  E_{2}^{t}=e_3 &
  E_{3}^{t}=e_4  &
E_{4}^{t}=te_1 &
E_{5}^{t}=e_5\\
\hline

\begin{array}{l|} \A_{23} \to \A_{24} \\
\hline\end{array}  &   
E_{1}^t= t e_1 + e_2&
E_{2}^t= 2 t e_3&
E_{3}^t= 2 t  e_2&
E_{4}^t= e_4 &
E_{5}^t= e_5 \\
\hline

\end{array}$$ }

{\tiny 
$$\begin{array}{|rcl|l|}
\hline
\multicolumn{4}{|c|}{\textrm{{\bf Table D. 
Non-degenerations of $5$-dimensional nilpotent associative commutative  algebras}}}  \\
\hline
\hline
\multicolumn{3}{|c|}{\textrm{Non-degeneration}} & \multicolumn{1}{|c|}{\textrm{Arguments}}\\
\hline
\hline

%\red{\C{02}} &\red{\not \to}& \red{\C{08}} & \red{\R=\{A_1A_3=0\}}\\ \hline

\A_{03} &\not \to & \A_{05} , \A_{07} & {\mathcal R}=
\left\{
 \begin{array}{l}    A_1^4=0, A_3^2 = 0,  A_1A_3 \subseteq A_5  \end{array}
\right\} \\ 
\hline

\A_{05  } &\not \to & \A_{15} & 
{\mathcal R}=
\left\{
 \begin{array}{l}   A_2A_3=0, c_{13}^4c_{22}^5=c_{13}^5c_{22}^4  \end{array}
\right\} \\ 
\hline

\A_{05}, \A_{06}, \A_{07} &\not \to &  \A_{21} &    {\mathcal R}=
\left\{
 \begin{array}{l}   % A_1A_4=0  
 \dim \operatorname{Ann}({\A_{21}}) = 1
 \end{array}
\right\} \\ 
\hline

\A_{07} &\not \to &  \A_{16},  \A_{18} &
{\mathcal R}=
\left\{
 \begin{array}{l}    A_1^2 \subseteq A_4, A_1^3=0  \end{array}
\right\} \\ 
\hline

\A_{10} &\not \to & \A_{19}, \A_{22} & 
{\mathcal R}=
\left\{
 \begin{array}{l}    A_2^2=0  \end{array}
\right\} \\ 
\hline

%\A_{11} &\not \to &  \A_{18}  & 
%{\mathcal R}=
%\left\{
 %\begin{array}{l}    A_1^3=0\end{array}
%\right\}  \\ 
%\hline

\A_{13} &\not \to &  \A_{12},  \A_{16} &    
{\mathcal R}=
\left\{
 \begin{array}{l}  A_1^2 \subseteq A_4,  A_1A_2 \subseteq A_5  \end{array}
\right\}, f_1=e_3, f_2=e_2, f_3=e_1, f_4=e_4, f_5=e_5 \\ 
\hline

%\A_{14} &\not \to & \A_{22} & 
%{\mathcal R}=
%\left\{
% \begin{array}{l}    A_2A_3=0  \end{array}
%\right\} \\ 
%\hline

\A_{16} &\not \to & \A_{12}, \A_{18}, \A_{22} & 
{\mathcal R}=
\left\{
 \begin{array}{l}    A_1A_3=0,  A_1^3=0 \end{array}
\right\} \\ 
\hline

\A_{21} &\not \to & \A_{20} & 
{\mathcal R}=
\left\{
 \begin{array}{l}    A_1^2 \subseteq A_5
 \end{array}
\right\} \\ 

\hline

\end{array}$$
}

\end{document}